\magnification 1200
\input amssym.def
\input amssym.tex
\parindent = 40 pt
\parskip = 12 pt

 at 22 true pt
\font \medheading =cmbx7 at 16 true pt
\font \small = cmr7 at 10  true pt
\font \heading = cmbx10 at 14 true pt
\font \smallheading = cmbx10 at 12 true pt

\def \R{{\bf R}}

\centerline{\medheading Fourier Transforms of Irregular Mixed Homogeneous }
\centerline{\medheading Hypersurface Measures }
\rm
\line{}
\line{}
\centerline{\heading Michael Greenblatt}
\line{}
\centerline {May 23, 2017}

\line{}
\line{}\baselineskip = 2 pt
\centerline {\bf Abstract} 
 {\narrower \noindent \small With the help of Van der Corput lemmas, decay estimates are proven for
 Fourier transforms of mixed homogeneous hypersurface measures with densities that can be quite irregular. The primary
results are local in nature, but can be extended to global theorems in an appropriate sense.
The estimates are sharp for a certain range of indices in the theorems.

\noindent MSC codes: 42B20 (primary), 42B10 (secondary). \hfill\break
 \par}

\line{}
\line{}
\line{} 
\noindent{\heading 1. Background and Theorem Statements}.\baselineskip = 12 pt

\vfootnote{}{This research was supported in part by NSF grant  DMS-1001070.}

\noindent We start with the following definition.

\noindent {\bf Definition 1.1.} Suppose $a_1,...,a_n$ are positive numbers. A function $f:\R^n \rightarrow \R$ is {\it mixed homogeneous of degree $(a_1,...,a_n)$} if for all $t > 0$ and all $(x_1,...,x_n) \in \R^n$ we have
$$f(t^{1 \over a_1}x_1,...,t^{1 \over a_n} x_n) = t f(x_1,...,x_n) \eqno (1.1)$$
Our first theorem will provide estimates on the decay rate of Fourier transforms of surfaces that are the graphs of mixed homogeneous functions. Specifically, we let $f(x)$ be a
function on $\R^n$ that is mixed homogeneous of degree $(a_1,...,a_n)$ and we let $g(x_1,...,x_n)$ be a locally integrable nonnegative function  that is mixed homogeneous of degree $(\rho a_1,..., \rho a_n)$ on $\R^n - \{0\}$ for some $\rho \in \R$
(which may be negative). We let $\phi(x_1,...,x_n)$ be an integrable function supported on the unit ball such that 
$$|\phi(x_1,...,x_n)| \leq g(x_1,...,x_n) \,\,\,\,\,\,\,\,\,{\rm and}\,\,\,\,\,\,\,\,\,\,\bigg|\sum_{i = 1}^n {x_i \over a_i}{\partial \phi \over \partial x_i} (x_1,...,x_n)\bigg| \leq g(x_1,...,x_n) \eqno (1.2)$$

A canonical example of a $\phi(x_1,...,x_n)$ satisfying $(1.2)$ would be $\phi(x_1,...,x_n) = \alpha(x_1,...,x_n)g(x_1,...,x_n)$, where $\alpha$ is smooth and supported in the unit ball. Note that the case where $\phi(x_1,...,x_n)$ is itself smooth is included in the situation where $\rho = 0$ and $g(x_1,...,x_n) = 1$.

\noindent We examine the Fourier transform of the surface measure of the graph of $f(x_1,...,x_n)$, weighted by $\phi(x_1,...,x_n)$, given by
$$T(\lambda_0,\lambda_1,...,\lambda_n) = \int_{\R^n} e^{i\lambda_0 f(x_1,...,x_n) + i\lambda_1 x_1 + ... + i\lambda_n x_n}\,\,\phi(x_1,...,x_n)
\,dx_1...\, dx_n \eqno (1.3)$$
Technically this is the Fourier transform of the surface measure at $(-\lambda_0,-\lambda_1,...,-\lambda_n)$, but to simplify notation we will consider $T(\lambda_0,\lambda_1,...,\lambda_n)$ as written here. Note that $f(x)$ and $g(x)$ (and therefore) $\phi (x)$ can be arbitrarily irregular in directions other than the mixed-radial direction; this is because our theorems (including the sharp estimates) will be proved with the use of Van der Corput lemmas in the
mixed-radial direction only. Hence $f(x)$, $g(x)$, and $\phi(x)$ can be badly behaved in other directions.

It should be pointed out that there is an extensive literature concerning the decay rate of Fourier transforms of hypersurface measures. When $\phi(x)$ is smooth and the Hessian of $f(x)$ has full rank everywhere, there are well-known asymptotic
formulas for the decay rate. We refer to chapter 8 of [S] for details. Whenever the
Hessian of the phase has positive rank, one can also obtain decay estimates generalizing the full rank case. We refer to [L]
as an example of this. It has been quite difficult to prove sharp estimates for general hypersurfaces (hence the restriction to mixed homogeneous surfaces here).  Many papers in this area 
have considered other classes of hypersurfaces such as the convex hypersurfaces considered in [BrNW] [CDMaM] [CMa] [NSeW] [Sc], 
or the case of scalar oscillatory integrals  (the situation where $\lambda_1 = .... = \lambda_n = 0$). We mention [ChKNo] [PY] [GPT] [Gr] [V] as examples of the latter.

Our first theorem 
regarding $T(\lambda_0,\lambda_1,...,\lambda_n)$ is as follows. Some cases of it follow from Theorem 14 of [ISa].

\noindent {\bf Theorem 1.1.} Let $f$ and $g$ be as above and let $v$ denote the number of distinct values in the set $\{a_1,...,a_n\}$. Assume that $a_i > 0$ and $a_i \neq 1$ for all $i$.  

\noindent {\bf a)} Let  $d\mu_g$ denote the measure $|g(x_1,...,x_n)|\, dx_1...\,dx_n$. Suppose 
 $\epsilon > 0$ is such that for some $B_{f, g, \epsilon} > 0$, for all $t > 0$ one has
$\mu_g (\{(x_1,...,x_n) \in [-1,1]^n: |f(x_1,...,x_n)| < t\}) \leq B_{f, g, \epsilon} t^{\epsilon}$.
 If $ \epsilon \neq {1 \over v + 1}$  there exists a $C_{f,g, \epsilon} > 0$ such that one has an 
estimate
$$|T(\lambda_0,\lambda_1,...,\lambda_n) | \leq C_{f,g, \epsilon}(1 + |\lambda_0|)^{-\min(\epsilon, {1 \over v + 1})} \eqno (1.4)$$
If $\epsilon ={1 \over v + 1}$,  then $(1.4)$ holds except with an additional factor of $\ln(2 + |\lambda_0|)$ on the right.  
Conversely, suppose $0 < \eta < 1$ is such that for all $\phi(x_1,...,x_n)$ supported in the unit ball satisfying $(1.2)$ one has 
an estimate
$$|T(\lambda_0,\lambda_1,...,\lambda_n) | \leq C_{f,g, \eta}(1 + |\lambda_0|)^{-\eta} \eqno (1.5)$$
Then for all $0 < \epsilon  < \eta$ one has an estimate of the form 
$$\mu_g (\{(x_1,...,x_n) \in [-1,1]^n: |f(x_1,...,x_n)| < t\}) \leq B_{f, g, \epsilon } t^{\epsilon} \eqno (1.6)$$
Hence if 
the supremum of all $\epsilon$ for which such a $B_{f, g,\epsilon}$ exists is at most ${1 \over v + 1}$, the exponent in
$(1.4)$ is sharp.

\noindent {\bf b)} Suppose $g(x_1,...,x_n) \in L^p([-1,1]^n)$, where $1 \leq p \leq \infty$.  Let $p'$ be the exponent 
conjugate to $p$; that is, ${1 \over p'} + {1 \over p} = 1$. If $p' \neq v +1 $ then for some $C_{f,g,p}' > 0$ one has an
estimate
$$|T(\lambda_0,\lambda_1,...,\lambda_n) | \leq C_{f,g,p}'(1 + |\lambda'|)^{-\min({1 \over p'}, {1 \over v + 1})} \eqno (1.7a)$$
Here $|\lambda'|$ denotes the magnitude of the vector $\lambda' = (\lambda_1,...,\lambda_n)$. If $p' = v + 1$ the same estimate holds with an additional factor of $(\ln (2 + |\lambda'|))^{1 \over v + 1}$ on the right-hand side.

\noindent {\bf c)} Suppose $g(x_1,...,x_n) \in  L^p([-1,1]^n)$ for some $1 \leq p \leq \infty$, and let $p'$ be the conjugate exponent as in part b). If $p' \neq v$ then for some $C_{f,g,p}'' > 0$ one has an 
estimate
$$|T(0,\lambda_1,...,\lambda_n) | \leq C_{f,g,p}''(1 + |\lambda'|)^{-\min({1 \over p'}, {1 \over v })} \eqno (1.7b)$$
If $p' = v$ the same estimate holds with an additional factor of $(\ln (2 + |\lambda'|))^{1 \over v}$.

Theorem 1.1 can be interpreted as follows. Suppose $a_i > 1$ for all $i$. The graph of $f(x_1,...,x_n)$ has a well-defined normal vector at the origin, given by $(0,...,0,1)$. Part a) of Theorem 1.1 gives estimates for the Fourier transform of the surface
 measure of this graph, weighted by 
$\phi(x_1,...,x_n)$, in terms of the component of $\lambda$ in this normal direction. The exponent is given by 
the growth rate in $t$ of this weighted measure of the portion of the surface with $|x_{n+1}| < t$.
Part b) gives estimates for tangential components of $\lambda$, while part c) estimates the Fourier transform of the weight function $\phi(x_1,...,x_n)$ itself similarly to part b).

When $v = 1$ (the homogeneous case) one can get an exponent as high as ${1 \over 2}$ in Theorem 1.1 a) and b). 
In this situation, with 
additional regularity assumptions on $f(x)$ and $g(x)$, one can use damping functions in conjunction with the methods of 
[SoS] to prove $L^p$ boundedness of  maximal averages associated with the hypersurfaces when $p > 2$. We refer to [ISa] 
for theorems of this kind. It might be interesting to explore what theorems of this nature are possible in the setting of this paper.

In the setting of Theorem 1.1a), one might ask if there are situations where $B_{f,g,\epsilon}$ exists for some $\epsilon \geq {1 \over v + 1}$ and where
the stronger estimate $|T(\lambda_0,\lambda_1,...,\lambda_n) | \leq C_{f,g, \epsilon}(1 + |\lambda_0|)^{-\epsilon}$ still holds. If one imposes 
extra regularity conditions on $f(x_1,...,x_n)$ and $\phi(x_1,...,x_n)$ in directions other than along the curves  $t \rightarrow (c_1t^a_1,...,c_nt^a_n)$ determined by the mixed homogeneity, 
it is not too hard
to find situations where this is the case. We refer to [ISa] for examples of this. But even in the absence of such regularity, it is conceivable that the worst
decay estimates one gets on these curves average out over $(c_1,...,c_n)$ so that one still obtains the 
stronger estimate. It would be interesting to try to determine in what situations this occurs, if any.

Although parts b) and c) of Theorem 1.1 are only sometimes sharp, when a homogeneity $a_i$ appears only once, one has
 sharp bounds on $|T(\lambda_0,\lambda_1,...,\lambda_n)|$ of the 
form $C|\lambda_i|^{-\epsilon}$ that are analogous to those given in Theorem 1.1 a). These bounds are given by the following
theorem; note that part c) gives a sharpness statement for both parts a) and b) .

\noindent {\bf Theorem 1.2.} Suppose the homogeneity $a_i$ is such that $a_j \neq a_i$ for all $j \neq i$. Let the measure $d\mu_g$ be as in Theorem 1.1. Suppose
$\epsilon_i > 0$ is such that for some $D_{f, g, \epsilon_i} > 0$, for all $t > 0$ one has
$\mu_g (\{(x_1,...,x_n) \in [-1,1]^n: |x_i| < t\}) \leq D_{f, g, \epsilon_i} t^{\epsilon_i}$.

 \noindent {\bf a)} If $ \epsilon_i \neq {1 \over v + 1}$  there exists a $E_{f,g, \epsilon_i} > 0$ such that one has an 
estimate
$$|T(\lambda_0,\lambda_1,...,\lambda_n) | \leq E_{f,g, \epsilon_i}(1 + |\lambda_i|)^{-\min(\epsilon_i, {1 \over v + 1})} \eqno (1.8a)$$
If $\epsilon_i ={1 \over v + 1}$,  then $(1.8a)$ holds except with an additional factor of $\ln(2 + |\lambda_i|)$ on the right.

\noindent {\bf b)} If  $ \epsilon_i \neq {1 \over v}$  there exists a $F_{f,g, \epsilon_i} > 0$ such that one has an 
estimate
$$|T(0,\lambda_1,...,\lambda_n) | \leq F_{f,g, \epsilon_i}(1 + |\lambda_i|)^{-\min(\epsilon_i, {1 \over v})} \eqno (1.8b)$$
If $\epsilon_i ={1 \over v}$,  then $(1.8b)$ holds except with an additional factor of $\ln(2 + |\lambda_i|)$ on the right.

\noindent {\bf c)} If $ 0 < \eta < 1$ is such that for all $\phi(x_1,...,x_n)$ satisfying $(1.2)$ one has an estimate
$$|T(0,\lambda_1,...,\lambda_n) | \leq F_{f,g, \eta}(1 + |\lambda_i|)^{-\eta} \eqno (1.9)$$
Then for all $0 < \epsilon  < \eta$ one has an estimate of the form 
$$\mu_g (\{(x_1,...,x_n) \in [-1,1]^n: |x_i| < t\}) \leq E_{f, g, \epsilon } t^{\epsilon}$$
---------------------------------------------------------------------------------------------------------------------

One can combine parts a) and b) of Theorem 1.1 to obtain decay estimates in terms of powers of $|\lambda|$, where $|\lambda|$ is the magnitude of the full vector $(\lambda_0,...,\lambda_n)$. For example, as an immediate consequence of Theorem 1.1 a) and b) we have the following.

\noindent {\bf Corollary 1.3.} Suppose  $\epsilon < {1 \over v +1} < {1 \over p'}$ and for some $B_{f, g, \epsilon} > 0$, for all $t > 0$ one has
$\mu_g (\{(x_1,...,x_n) \in [-1,1]^n: |f(x_1,...,x_n)| < t\}) \leq B_{f, g, \epsilon} t^{\epsilon}$. Then 
for some $C_{f,g,\epsilon} > 0$ we have an estimate
$$|T(\lambda_0,\lambda_1,...,\lambda_n) | \leq C_{f,g,\epsilon} (1 + |\lambda|)^{-\epsilon} \eqno (1.10) $$
If the supremum of the $\epsilon$ for which such a $B_{f, g, \epsilon}$ exists is given by $\delta \leq {1 \over v + 1}$, then the
exponent in the right-hand side of $(1.10)$ cannot be taken to be greater than $\delta$.

Next,
we will have a global analogue of the local result of Corollary 1.3 for smooth $\phi(x_1,...,x_n)$. To ensure that all
integrals are well-defined, our global analogue is as follows. While $\phi(x_1,...,x_n)$ will still be compactly
supported, its support can be arbitrarily large, and we will obtain uniform bounds on  $|T(\lambda_0,\lambda_1,...,\lambda_n) |$ over all $\phi(x_1,...,x_n)$ that satisfy the following decay conditions. 
 We stipulate that for some $k > \sum_{i=1}^n {1 \over a_i} - \epsilon$ and some $A > 0$ we have 
$$|\phi(x_1,...,x_n)| \leq A(1 + \sum_{i=1}^n |x_i|^{a_i})^{-k}  ,\,\,\,\,\,\,\,\,\,\,\,\,\,\,\, \bigg|\sum_{i = 1}^n {x_i\over a_i} {\partial \phi \over \partial x_i} (x_1,...,x_n)\bigg| \leq  A(1 + \sum_{i=1}^n |x_i|^{a_i})^{-k} \eqno (1.11)$$

\noindent {\bf Theorem 1.4.}
Suppose $(1.11)$ holds and  suppose $0 < \epsilon < {1 \over v + 1}$ is such that for some constant $C_{f,\epsilon}$ the Lebesgue measure of the set $\{(x_1,...,x_n) \in [-1,1]^n: |f(x_1,...,x_n)| < t\}$ is bounded by $C_{f,\epsilon} t^{\epsilon}$ for all $t > 0$. Then 
for some $D_{f,A,\epsilon} >0 $ we have the estimate
$$|T(\lambda_0,\lambda_1,...,\lambda_n) | \leq D_{f,A,\epsilon} (1 + |\lambda|)^{-\epsilon} \eqno (1.12)$$
The exponent in Theorem 1.4  is best possible when  $|\{x \in [-1,1]^n: |f(x)| < t\}| \sim t^{\epsilon}$ for small $t$; if $\phi$ is nonnegative  with $\phi(0,...,0) \neq  0$ and is supported on a small enough neighborhood of the origin, then the proof of the sharpness statement  of Theorem 1.1a) at the end of this section shows that the exponent $\epsilon$ cannot be improved for this $\phi(x)$.

\noindent {\smallheading Examples.}

\noindent {\bf Example 1.} Suppose $f(x_1,...,x_n)$  is a bounded function on $[-1,1]^n$ such that for some $C, r > 0$, for each 
$t > 0$ one has
 $|\{(x_1,...,x_n) \in [-1,1]^n: |f(x_1,...,x_n)| < t\}| < C t^r$ (absolute values here denotes Lebesgue measure).
Let $g(x_1,...,x_n) = |f(x_1,...,x_n)|^{\rho}$ for some $\rho > -r$, and let $\phi(x_1,...,x_n) = \alpha(x_1,...,x_n) g(x_1,...,x_n) =  \alpha(x_1,...,x_n)|f(x_1,...,x_n)|^{\rho}$, where $\alpha(x_1,...,x_n)$ is a smooth compactly supported function on $(-1,1)^n$. Note that geometrically, a damping function  $|f(x_1,...,x_n)|^{\rho}$ represents $|x_{n+1}|^{\rho}$ when the surface is viewed in $\R^{n+1}$.

\noindent Then we have that
$\mu_g (\{(x_1,...,x_n) \in [-1,1]^n: |f(x_1,...,x_n)| < t\}) $ is given by 
$$ \int_{\{(x_1,...,x_n)\in [-1,1]^n: |f(x_1,...,x_n)| < t\}} |f(x_1,...,x_n)|^{\rho}\,dx_1...\,dx_n \eqno (1.13)$$
$$= \sum_{k = 0}^{\infty}\int_{\{(x_1,...,x_n)\in [-1,1]^n: \,2^{-k-1} t \,\leq\,|f(x_1,...,x_n)| < 2^{-k} t\}} |f(x_1,...,x_n)|^{\rho}\,dx_1...\,dx_n $$
$$\leq \sum_{k = 0}^{\infty} C(2^{-k}t)^r(2^{-k}t)^{\rho} \eqno (1.14)$$
$$ < C't^{\rho + r} \eqno (1.15)$$
Hence Theorem 1.1 part a) holds with $\epsilon = \rho + r$. As for parts b) and c), if $\rho \geq  0$, then one can take $p = \infty$.
If $\rho < 0$, then one can take $p$ to be any number less than $ {r \over  |\rho|}$.

 \noindent As for Theorem 1.2, note that one has that 
$\mu_g (\{(x_1,...,x_n) \in [-1,1]^n: |x_i| < t\}) $ is given by 
$$\int_{\{(x_1,...,x_n)\in [-1,1]^n: |x_i |< t\}} |f(x_1,...,x_n)|^{\rho}\,dx_1...\,dx_n \eqno (1.16)$$
Thus the exponents in Theorem 1.2 are determined by the growth rate in $t$ of the integral $(1.16)$.

\noindent {\bf Example 2.} Suppose we are in the homogeneous case; that is, there is an $a > 0$, $a \neq 1$, such that $a_i = a$ for all $i$. 
Then $v = 1$, which gives the best possible exponents in Theorem 1.1; one can get a decay rate up to $(1 + |\lambda_0|)^{-{1 \over 2}}$ and  $(1 + |\lambda'|)^{-{1 \over 2}}$ in parts a) and b) respectively. Suppose $g(x_1,...,x_n) = |x_1^{\alpha_1}...
x_n^{\alpha_n}|$ for some nonnegative integers $\alpha_i$, so that  for example $\phi(x_1,...,x_n)$ can be of the form $\alpha(x_1,...,x_n) x_1^{\alpha_1}...x_n^{\alpha_n}$ where $\alpha(x_1,...,x_n)$ is a smooth compactly supported function on $(-1,1)^n$. 

\noindent  Then $\mu_g (\{(x_1,...,x_n) \in [-1,1]^n: |f(x_1,...,x_n)| < t\}) $  is equal to
$$ \int_{\{(x_1,...,x_n)\in [-1,1]^n: \,|f(x_1,...,x_n)| < t\}} |x_1^{\alpha_1}...x_n^{\alpha_n}|\,dx_1...\,dx_n \eqno (1.17)$$
If each $\alpha_i$ is even for example, we can change variables to $y_i = x_i^{\alpha_i + 1}$ for all $i$. Then the above becomes
$$\prod_{i = 1}^n {1 \over \alpha_i + 1} \int_{\{(y_1,...,y_n) \in [-1,1]^n:\, |f(y_1^{1 \over \alpha_1 + 1},...,y_n^{1 \over \alpha_n + 1})| < t\}}\,dy_1...\,dy_n $$
Thus $\mu_g (\{(x_1,...,x_n) \in [-1,1]^n: |f(x_1,...,x_n)| < t\}) $ is expressible in terms of the distribution function of $|f(x_1^{1 \over \alpha_1 + 1},...,x_n^{1 \over \alpha_n + 1})|$. There is an analogous interpretation when the $\alpha_i$ are allowed to be odd. Note that the exponent $p$ in parts b) and c) of Theorem 1.1 in this example can be taken to be $\infty$.

\noindent {\bf Example 3.} Suppose we are in the opposite situation from Example 2, and each $a_i$ is different, so that $v = n$.
Then for a given $i$, Theorem 1.2 provides estimates of the form $|T(\lambda)| \leq F_{f,g, \epsilon_i}(1 + |\lambda_i|)^{-\min(\epsilon_i, {1 \over n + 1})} $ (when a positive $\epsilon_i$ exists) and one has an overall estimate $|T(\lambda)| \leq F_{f,g, \epsilon}(1 + |\lambda|)^{-\epsilon}$ where $\epsilon = \min(\epsilon_1,...,\epsilon_n, {1 \over n + 1})$. If for some $i$ the exponent $\epsilon_i$ cannot be taken to be greater than ${1 \over n+1}$, then by Theorem 1.2c) the minimum over $i$ of the supremal $\epsilon_i$ gives the supremum of the exponents $\eta$ for 
which we have an estimate of the form $|T(\lambda)| \leq F_{f,g, \eta}(1 + |\lambda|)^{-\eta}$.

\noindent {\smallheading Proof of sharpness statements in Theorem 1.1 and 1.2.}

The sharpness statements in Theorem 1.1a) and Theorem 1.2c) can be proved in relatively short order. We focus our attention on
the situation in Theorem 1.1a). Suppose $0 < \eta < 1$ is such that $(1.5)$ holds for all $\phi(x_1,...,x_n)$ 
supported in the unit ball satisfying $(1.2)$. Let $\alpha(x)$ be a bump function on $\R$  whose Fourier transform is nonnegative, compactly supported, and equal to 1 on a 
neighborhood of the origin, and let $N$ be a large positive number. If $0 <  \epsilon < \eta$, then $(1.5)$ implies that for 
some constant $A_{f,g,\epsilon}$ one has
$$\int_{\R}|T(\lambda_0,0,...,0)| |\lambda_0|^{\epsilon  - 1}\alpha(N \lambda_0)\,d \lambda_0 < A_{f,g,\epsilon} \eqno (1.18)$$
In view of the definition $(1.3)$ for $T(\lambda_0,...,\lambda_n)$, this implies that
$$\bigg|\int_{\R^{n+1}} e^{i\lambda_0 f(x_1,...,x_n)}\,\phi(x_1,...,x_n) |\lambda_0|^{\epsilon  - 1}\alpha(N \lambda_0)\,d\lambda_0 \,dx_1...\, dx_n\bigg|  <  A_{f,g,\epsilon} \eqno (1.19)$$
If we do the integral in $\lambda_0$ first in $(1.19)$, we get
$$\bigg|\int_{\R^n}\beta_N (f(x_1,...,x_n)) \,\phi(x_1,...,x_n)\,dx_1...\, dx_n\bigg| < A_{f,g,\epsilon}' \eqno (1.20)$$
Here $\beta_N(y)$ is the convolution of $|y|^{-\epsilon}$ with ${1 \over N} \hat{\alpha}({y \over N})$.
Taking $\phi(x)= \gamma(x)g(x)$ for some nonnegative bump function $\gamma(x)$ equal to one on a neighborhood of the origin we get
$$\bigg|\int_{\R^n}\beta_N (f(x))\gamma(x)g(x)\,dx_1...\, dx_n\bigg| < A_{f,g,\epsilon}''  \eqno (1.21)$$
Letting $N \rightarrow \infty$ gives
$$\int_{\R^n}|f(x)|^{-\epsilon}\gamma(x)g(x)\,dx_1...\, dx_n < A_{f,g,\epsilon}'''  \eqno (1.22)$$
Since $\gamma(x)$ is equal to 1 on a neighborhood of the origin, mixed homogeneity gives
$$\int_{[-1,1]^n}|f(x)|^{-\epsilon} g(x)\,dx_1...\, dx_n < \infty  \eqno (1.23)$$
In other words, $|f(x)|^{-\epsilon}$ is in $L^1([-1,1]^n)$ with respect to the measure $d\mu_g$. Hence it is in weak $L^1$, 
and we have the existence of a constant $G_{f,g,\epsilon}$ such that
$$\mu_g (\{(x_1,...,x_n) \in [-1,1]^n: |f(x_1,...,x_n)|^{-\epsilon} > t\}) \leq G_{f, g, \epsilon } {1 \over t} \eqno (1.24)$$
Replacing $t$ by $t^{-\epsilon}$, $(1.24)$ gives the sharpness statement $(1.6)$ as needed. This gives the sharpness statement
in Theorem 1.1a).

Theorem 1.2c) is proved in exactly the same way, replacing  the $\lambda_0$ variable by the $\lambda_i$ variable and
the function $f(x_1,...,x_n)$ by $x_i$.

\noindent {\heading 2. Proofs of Theorems 1.1, 1.2, and 1.4.}

\noindent We start with the well-known Van der Corput lemma (see p. 334 of [S]).

\noindent {\bf Lemma 2.1.} Suppose $h(x)$ is a $C^k$ function on the interval $[a,b]$ with $|h^{(k)}(x)| > A$ on $[a,b]$ for
some $A > 0$. Let $\phi(x)$ be $C^1$ on $[a,b]$. If $k \geq 2$ there is a constant $c_k$ depending only on $k$ such that
$$\bigg|\int_a^b e^{ih(x)}\phi(x)\,dx\bigg| \leq c_kA^{-{1 \over k}}\bigg(|\phi(b)| + \int_a^b |\phi'(x)|\,dx\bigg) \eqno (2.1)$$
If $k =1$, the same is true if we also assume that $h(x)$ is $C^2$ and $h'(x)$ is monotone on $[a,b]$. 

\noindent  We will also make use of the following variant of Lemma 2.1 for $k = 1$.

\noindent {\bf Lemma 2.2.} Suppose the hypotheses of Lemma 2.1 hold with $k = 1$, except  instead of assuming that $h'(x)$ is monotone on $[a,b]$ we assume that $|h''(x)|
< {B  \over(b-a)}A$ for some constant $B > 0$. Then we have
$$\bigg|\int_a^b e^{ih(x)}\phi(x)\,dx\bigg| \leq  A^{-1}\bigg(\int_a^b |\phi'(x)|\,dx+ (B+2) 
\sup_{[a,b]}|\phi(x)| \bigg) \eqno (2.2)$$

\noindent {\bf Proof.} We write $e^{ih(x)} = (h'(x)e^{ih(x)}) {1 \over h'(x)}$ and integrate by parts in the integral being estimated,
integrating the $h'(x)e^{ih(x)}$ factor to $e^{ih(x)}$ and differentiating ${\phi(x) \over h'(x)}$. We obtain
$$\int_a^b e^{ih(x)}\phi(x)\,dx =  e^{ih(b)}{\phi(b) \over h'(b)} - e^{ih(a)}{\phi(a) \over h'(a)} - \int_a^b e^{ih(x)}{d \over dx}\bigg({\phi(x) \over h'(x)}\bigg)\,dx$$
$$= e^{ih(b)}{\phi(b) \over h'(b)} - e^{ih(a)}{\phi(a) \over h'(a)}- \int_a^b e^{ih(x)}{\phi'(x) \over h'(x)}\,dx +  \int_a^b e^{ih(x)}{\phi(x) h''(x)\over (h'(x))^2}\,dx \eqno (2.3)$$
The condition that $|h'(x)| > A$ ensures that each of the two boundary terms is bounded in absolute value by $A^{-1} 
\sup_{[a,b]}|\phi(x)|$. As for the first integral term, taking absolute values of the integrand and inserting $|h'(x)| > A$ gives that
this term is bounded in absolute value by $A^{-1}\int_a^b |\phi'(x)|\,dx$. In the second integral term, we use that $|h''(x)| < 
{B  \over(b-a)}A$ and $|h'(x)| > A$, resulting in 
$$\bigg|{\phi(x) h''(x)\over (h'(x))^2}\bigg| \leq B{A^{-1}\over (b-a)}|\phi(x)|$$
$$ \leq B{A^{-1} \over (b-a)}\sup_{[a,b]}|\phi(x)| \eqno (2.4)$$
Thus the second integral term is bounded by $BA^{-1}\sup_{[a,b]}|\phi(x)|$. Adding the bounds for the different terms gives
us the bounds on the right-hand side of $(2.2)$ and we are done.

\noindent {\bf Beginning of the Proof of Theorem 1.1.} 

\noindent We will prove the bounds of Theorem 1.1 for $|T'(\lambda_0,\lambda_1,...,\lambda_n)| =$
$$\int_{[0,1]^n} e^{i\lambda_0 f(x_1,...,x_n) + i\lambda_1 x_1 + ... + i\lambda_n x_n}\phi(x_1,...,x_n)
\,dx_1...\, dx_n \eqno (2.5)$$
The analogous bounds for the other octants of the integral of $(1.3)$ follow from changing variables to $(\pm x_1,...,\pm x_n)$
and then using the same argument. We change variables in $(2.5)$ via $(x_1,...,x_n) = (y_1^{1 \over a_1},...,y_n^{1 \over a_n})$,
obtaining
$$\int_{[0,1]^n} e^{i\lambda_0 f(y_1^{1 \over a_1},...,y_n^{1 \over a_n}) + i\lambda_1 y_1^{1 \over a_1} + ... + i\lambda_n y_n^{1 \over a_n}}$$
$$\times {1 \over a_1...a_n}y_1^{{1 \over a_1}- 1}...y_n^{{1 \over a_n}- 1}\phi(y_1^{1 \over a_1},...,y_n^{1 \over a_n})
\,dy_1...\, dy_n \eqno (2.6)$$
We define
$$F(y_1,...,y_n) = f(y_1^{1 \over a_1},...,y_n^{1 \over a_n})$$
$$\psi(y_1,...,y_n) = {1 \over a_1...a_n}y_1^{{1 \over a_1}- 1}...y_n^{{1 \over a_n}- 1}\phi(y_1^{1 \over a_1},...,y_n^{1 \over a_n})$$
$$\bar{g}(y_1,...,y_n) = {1 \over a_1...a_n}y_1^{{1 \over a_1}- 1}...y_n^{{1 \over a_n}- 1}g(y_1^{1 \over a_1},...,y_n^{1 \over a_n}) \eqno (2.7)$$
Here $g(x_1,...,x_n)$ is as in $(1.2)$. 
Then  $F$ is homogeneous of degree one, $\bar{g}$ is homogeneous of degree $\sum_{i = 1}^n({1 \over a_i} -1 )  + \rho$, where $\rho$ is as in the definition of $g(x_1,...,x_n)$, and a direct calculation reveals $\psi(y_1,...,y_n)$
 satisfies $(1.2)$ with each $a_i = 1$ and where $g(x_1,...,x_n)$ replaced by a constant multiple of $\bar{g}(y_1,...,y_n)$ which we henceforth refer to as $\tilde{g}(y_1,...,y_n)$.This constant depends only on the $a_i$. We next
rewrite $(2.6)$ as
$$\int_{[0,1]^n} e^{i\lambda_0 F(y_1,...,y_n) + i\lambda_1 y_1^{1 \over a_1} + ... + i\lambda_n y_n^{1 \over a_n}}\psi(y_1,...,y_n)
\,dy_1...\, dy_n \eqno (2.8)$$
Let $S_n^+$ denote $\{(v_1,...,v_n) \in S^n: \,v_i >\,0{\rm\,\,for\,\,all\,\,}i\}$, where $S^n$ denotes the unit sphere. In polar coordinates, $(2.6)$ becomes a constant depending only on $n$ times
$$\int_{S_n^+} \int_0^{\sqrt{n}} e^{i\lambda_0 F(rv_1,...,rv_n) + i\lambda_1 (rv_1)^{1 \over a_1} + ... + i\lambda_n (rv_n)^{1 \over a_n}}r^{n-1} \psi(rv_1,...,rv_n)
\,dr\,dv_1...\, dv_n \eqno (2.9)$$
The inner integral here goes to $\sqrt{n}$ due to the integrand being supported on $[0,1]^n$. Note that the condition $(1.2)$
when  $a_i = 1$ for all $i$ translates into
$$|\psi(rv_1,...,rv_n)| \leq \tilde{g}(rv_1,...,rv_n)\,\,\,\,\,\,\,\,\,\,\,\,\,\,\,\,\,\,|\partial_r \psi(rv_1,...,rv_n)| \leq  {1 \over r}
\tilde{g}(ry_1,...,ry_n) \eqno (2.10a)$$
Since $\tilde{g}$ is homogeneous of degree $t = \sum_{i = 1}^n({1 \over a_i} -1 )  + \rho$, $(2.10a)$ can further be rewritten as 
$$|\psi(rv_1,...,rv_n)| \leq  \tilde{g}(v_1,...,v_n)r^t\,\,\,\,\,\,\,\,\,\,\,\,\,\,\,\,\,\,|\partial_r \psi(rv_1,...,rv_n)| \leq \tilde{g}(v_1,...,v_n)r^{t -1  }\eqno (2.10b)$$
Next, we enumerate the $v$ elements of the set $\{{1 \over a_1},...,{1 \over a_n}\}$ as $\{b_1,...,b_v\}$, and we denote the
coefficient $\sum_{\{k: {1 \over a_k} = b_i\}} \lambda_kv_k^{1 \over a_k}$ of $r^{b_i}$ in the phase function of $(2.9)$ by $\mu_i$. Using also that $F$ is homogeneous of degree one, $(2.9)$ can therefore be rewritten as 
$$\int_{S_n^+} \int_0^{\sqrt{n}} e^{i\lambda_0 F(v_1,...,v_n)r + i\mu_1r^{b_1}  + ... + i\mu_vr^{b_v}}r^{n-1} \psi(rv_1,...,rv_n)
\,dr\,dv_1...\, dv_n \eqno (2.11)$$
Since Theorem 1.1 stipulates that no $a_i = 1$, we also have that no $b_i = 1$. We next divide the inner integral of $(2.11)$ dyadically, writing it as
$$\sum_{m\geq \log(-\sqrt{n})} \int_{2^{-m - 1}}^{2^{-m}} e^{i\lambda_0 F(v_1,...,v_n)r + i\mu_1r^{b_1}  + ... + i\mu_vr^{b_v}}r^{n-1} \psi(rv_1,...,rv_n)
\,dr \eqno (2.12)$$
Let $P(r)$ denote the phase function $\lambda_0 F(v_1,...,v_n)r + \mu_1r^{b_1}  + ... +  \mu_vr^{b_v}$ in $(2.12)$. In the following, let $b_0 = 1$ so that
the power of $r$ in the first term of $P(r)$ is on the same footing as the power of $r$ in the other terms of $P(r)$. Let ${\bf z}$
denote the $v + 1$ by $1$ column vector whose $i$th entry is $r^i \partial_r^i P(r)$. Let ${\bf M}$ denote the $v + 1$ by $v + 1$ matrix whose $ij$
entry is $b_{j-1}(b_{j-1}- 1)...(b_{j-1}- i+ 1)$, and let ${\bf w}$ denote the $v + 1$ by $1$ column vector whose first entry is $\lambda_0 F(v_1,...,v_n)r$ and whose $i$th entry for $i > 1$ is $\mu_{i-1} r^{b_{i-1}}$. Then we have the identity
$${\bf z} = {\bf M}{\bf w} \eqno (2.13)$$
The matrix ${\bf M}$ is invertible since after elementary row operations it can be converted into the Vandermonde matrix whose 
$i$th row is given by $(b_0^i,...,b_n^i)$, and all the $b_j$ are distinct. As a result, there is a constant $c$ depending only on 
the $b_i$ such that there is always some $i$ for which 
$$|z_i| \geq c|{\bf w}| \eqno (2.14a)$$
More explicitly, there is a constant $c'$ depending only on the $b_j$ such that for any $r$ there is some $i$ with $1 \leq i \leq v + 1$ such that
$$|\partial_r^i (\lambda_0 F(v_1,...,v_n)r + \mu_1r^{b_1}  + ... + \mu_vr^{b_v})| \geq  c'{1 \over r^i}\bigg( |\lambda_0 F(v_1,...,v_n)r|
+ \sum_{l = 1}^v |\mu_lr^{b_l}|\bigg) \eqno (2.14b)$$
Furthermore, by directly bounding each term, one immediately has that there is a constant $c''$ depending on the $b_j$ such that 
$$|\partial_r^{i + 1}(\lambda_0 F(v_1,...,v_n)r + \mu_1r^{b_1}  + ... + \mu_vr^{b_v})| \leq  c''{1 \over r^{i + 1}}\bigg( |\lambda_0 F(v_1,...,v_n)r| + \sum_{l = 1}^v |\mu_lr^{b_l}|\bigg) \eqno (2.15)$$
Hence $(2.14b)$ implies that there is some $c_1 > 0$ depending on the $b_j$ such that for each given $r_0$, for $r \in [(1 - c_1)r_0, (1 + c_1)r_0]$ one has
$$|\partial_r^i (\lambda_0 F(v_1,...,v_n)r + \mu_1r^{b_1}  + ... + \mu_vr^{b_v})| \geq  {c' \over 2}{1 \over r^i}\bigg( |\lambda_0 F(v_1,...,v_n)r_0| + \sum_{l= 1}^v |\mu_lr_0^{b_l}|\bigg) \eqno (2.16)$$
Going back to $(2.12)$, we see that $(2.16)$ implies that a given term of $(2.12)$ may be written as the sum of at most ${1 \over c_1}$ terms on which $(2.16)$ holds for a single $i$. To this end, we write 
$$\int_{2^{-m - 1}}^{2^{-m}} e^{i\lambda_0 F(v_1,...,v_n)r + i\mu_1r^{b_1}  + ... + i\mu_vr^{b_v}}r^{n-1} \psi(rv_1,...,rv_n)
\,dr$$
$$  = \sum_{k = 0}^{1 \over c_1} \int_{I_k} e^{i\lambda_0 F(v_1,...,v_n)r + i\mu_1r^{b_1}  + ... + i\mu_vr^{b_v}}r^{n-1} \psi(rv_1,...,rv_n)
\,dr  \eqno (2.17)$$
Here the $I_k$ denote intervals on which $(2.16)$ holds for a single $i$ on the interval $I_k$. Next, we use $(2.16)$ to apply the 
Van der Corput lemma on each term of $(2.17)$. We apply Lemma 2.1 if $i > 1$ and Lemma 2.2 if $i = 1$, using $(2.10b)$ to bound $|\psi(rv_1,...,rv_n)|$ and its $r$ derivative.
The result is
$$ \bigg| \int_{I_k} e^{i\lambda_0 F(v_1,...,v_n)r + i\mu_1r^{b_1}  + ... + i\mu_vr^{b_v}}r^{n-1} \psi(rv_1,...,rv_n)\,dr \bigg|$$
$$\leq Cr_0 \bigg(|\lambda_0 F(v_1,...,v_n)r_0| + \sum_{l= 1}^v |\mu_lr_0^{b_l}|\bigg)^{-{1 \over i}} r_0^{n-1}|\tilde{g}(v_1,...,v_n)| r_0^{t + 1}\eqno (2.18)$$
Here $C$ depends on the $b_j$ and $v$. It should be pointed out that the bound $(2.18)$ for a term of $(2.12)$ also follows from
 the corollary to the proposition of section 3 of [RS], but we include the somewhat related argument here for completeness. 

\noindent Since $I_k$ is ${1 \over c_1}$ of a dyadic interval, the right-hand side of $(2.18)$ is in turn bounded by 
$$C' \int_{I_k}   \bigg(|\lambda_0 F(v_1,...,v_n)r| + \sum_{l = 1}^v |\mu_lr^{b_l}|\bigg)^{-{1 \over i}} r^{n-1}|\tilde{g}(v_1,...,v_n)| r^t \,dr \eqno (2.19)$$
On the other hand, simply by taking absolute values of the integrand and integrating, in view of $(2.10b)$ one has
$$\int_{2^{-m - 1}}^{2^{-m}} e^{i\lambda_0 F(v_1,...,v_n)r + i\mu_1r^{b_1}  + ... + i\mu_vr^{b_v}}r^{n-1} \psi(rv_1,...,rv_n)
\,dr$$
$$\leq  C''  \int_{2^{-m - 1}}^{2^{-m}}  r^{n-1}|\tilde{g}(v_1,...,v_n)| r^t \,dr \eqno (2.20)$$
Hence one may combine $(2.19)$ and $(2.20)$ to conclude that
$$\int_{2^{-m - 1}}^{2^{-m}} e^{i\lambda_0 F(v_1,...,v_n)r + i\mu_1r^{b_1}  + ... + i\mu_vr^{b_v}}r^{n-1} \psi(rv_1,...,rv_n)
\,dr $$
$$\leq C_3\int_{2^{-m - 1}}^{2^{-m}}\min\bigg(1,  \bigg(|\lambda_0 F(v_1,...,v_n)r| + \sum_{l = 1}^v |\mu_lr^{b_l}|\bigg)^{-{1 \over i}}\bigg)  r^{n-1}|\tilde{g}(v_1,...,v_n)| r^t \,dr \eqno (2.21)$$
This is maximized when $i$ is as large as possible, namely $i = v + 1$. Hence $(2.21)$ is bounded by
$$\leq C_3\int_{2^{-m - 1}}^{2^{-m}}\min\bigg(1,  \bigg(|\lambda_0 F(v_1,...,v_n)r| + \sum_{l = 1}^v |\mu_lr^{b_l}|\bigg)^{-{1 \over v + 1}}\bigg)  r^{n-1}|\tilde{g}(v_1,...,v_n)| r^t \,dr \eqno (2.22)$$
Adding over all $k$, we get 
$$\int_{2^{-m - 1}}^{2^{-m}} e^{i\lambda_0 F(v_1,...,v_n)r + i\mu_1r^{b_1}  + ... + i\mu_vr^{b_v}}r^{n-1} \psi(rv_1,...,rv_n)
\,dr$$
$$\leq C' \int_{2^{-m - 1}}^{2^{-m}}   \bigg(|\lambda_0 F(v_1,...,v_n)r| + \sum_{l = 1}^v |\mu_lr^{b_l}|\bigg)^{-{1 \over v + 1}} r^{n-1}|\tilde{g}(v_1,...,v_n)| r^t \,dr \eqno (2.23)$$
Next, we add $(2.23)$ over all $m$ and insert it back into $(2.11)$. We obtain 
$$\bigg|\int_{S_n^+} \int_{r \geq 0}e^{i\lambda_0 F(v_1,...,v_n)r + i\mu_1r^{b_1}  + ... + i\mu_vr^{b_v}}r^{n-1} \psi(rv_1,...,rv_n)
\,dr\,dv_1...\, dv_n\bigg| $$
$$\leq C_4 \int_{(rv_1,...,rv_n) \in [0,1]^n}  \min\bigg(1,  \bigg(|\lambda_0 F(v_1,...,v_n)r| + \sum_{l = 1}^v |\mu_lr^{b_l}|\bigg)^{-{1 \over v + 1}}\bigg)$$
$$\times  r^{n-1}|\tilde{g}(v_1,...,v_n)| r^t \,dr \, dv_1...\,dv_n\eqno (2.24)$$
One technical point in the right-hand side of $(2.24)$ is that for a given $(v_1,...,v_n)$ the maximum possible $r$ might be in the middle of a dyadic interval and thus adding over all $m$ and integrating does not exactly give $(2.24)$. However, since the integrand on varies by at most a fixed factor over any such dyadic interval the estimate $(2.24)$ will still hold.

\noindent We next convert back from polar into rectangular coordinates in the right-hand side of $(2.24)$, obtaining 
$$\bigg|\int_{[0,1]^n} e^{i\lambda_0 F(y_1,...,y_n) + i\lambda_1 y_1^{1 \over a_1} + ... + i\lambda_n y_n^{1 \over a_n}}\psi(y_1,...,y_n)
\,dy_1...\, dy_n \bigg| $$
$$\leq C_5\int_{[0,1]^n} \min\bigg(1,  \bigg(|\lambda_0 F(y_1,...,y_n)| + \sum_{l = 1}^v \big|\sum_{\{m:\, b_m= b_l\}} \lambda_m y_m^{b_l}\big|\bigg)^{-{1 \over v + 1}}\bigg)$$
$$\times   |\tilde{g}(y_1,...,y_n)|\, dy_1...\,dy_n\eqno (2.25)$$
Converting this back into the original $x$-coordinates of $(2.5)$ gives
$$\bigg| \int_{[0,1]^n}  e^{i\lambda_0 f(x_1,...,x_n) + i\lambda_1 x_1+ ... + i\lambda_n x_n}\phi(x_1,...,x_n)
\,dx_1...\, dx_n\bigg| $$
$$\leq C_6\int_{[0,1]^n} \min\bigg(1,  \bigg(|\lambda_0 f(x_1,...,x_n)| + \sum_{l = 1}^v \big|\sum_{\{m:\, a_m = a_l\}} \lambda_m x_m\big|\bigg)^{-{1 \over v + 1}}\bigg)$$
$$\times   |g(x_1,...,x_n)|\, dx_1...\,dx_n\eqno (2.26)$$
We are now are in a position to prove Theorem 1.1. For part a), we will use $|\lambda_0 f(x_1,...,x_n)| $ in the minimum in
$(2.26)$ and for part b) we will use the $|\sum_{\{m:\, a_m = a_i\}} \lambda_m x_m\big|$ terms. Part c) is essentially a repeat
of part b), except that with the absence of a $\lambda_0 f(x_1,...,x_n)$ term in the phase we will be able to replace $v + 1$ by $v$. We start with part a).

\line{}

\noindent {\bf Proof of part a) of Theorem 1.1.} 

From now until the  end of the proof of part a) of Theorem 1.1 we use the notation $C$ to denote a constant that may depend on $f, g, $ and $\epsilon$ as in the statement of Theorem 1.1a). Let $d\mu_g$ be the measure $|g(x_1,...,x_n)|\, dx_1...\,dx_n$ as in the
statement of Theorem 1.1. Letting
$T'(\lambda_0,...,\lambda_n)$ denote the integral in the left-hand side $(2.26)$ as before, by $(2.26)$ we have that
$$|T'(\lambda_0,...,\lambda_n)| \leq C\int_{[0,1]^n}\min(1, |\lambda_0 f(x_1,...,x_n)|^{-{1 \over v + 1}})\,d\mu_g \eqno (2.27)$$
$$= \mu_g(\{(x_1,...,x_n) \in [0,1]^n : |f(x_1,...,x_n)| < |\lambda_0|^{-1} \}) $$
$$+ |\lambda_0|^{-{1 \over v + 1}} \int_{ \{(x_1,...,x_n) \in [0,1]^n : |f(x_1,...,x_n)| > {1 \over |\lambda_0|}\}} |f(x_1,...,x_n)|^{-{1 \over v + 1}}d\mu_g  \eqno (2.28)$$
By the characterization of integrals of powers of functions in terms of their distribution functions, applied to $|f(x_1,...,x_n)|^{-1}$,
the integral in $(2.28)$ is equal to
$${1 \over v + 1}\int_{|\lambda_0|^{-1}}^{\infty}t^{-{{1 \over v + 1} - 1}}\mu_g(\{(x_1,...,x_n) \in [0,1]^n: |\lambda_0|^{-1}
< |f(x_1,...,x_n)| < t\})\,dt \eqno (2.29)$$
Thus if $\epsilon$ satisfies an estimate $\mu_g (\{(x_1,...,x_n) \in [-1,1]^n: |f(x_1,...,x_n)| < t\}) \leq B_{f, g, \epsilon} t^{\epsilon}$, then  $(2.29)$ is bounded by
$$C \int_{|\lambda_0|^{-1}}^{\infty}t^{-{{1 \over v + 1} - 1}}\min(1, t^{\epsilon}) \,dt \eqno (2.30)$$
We can put the minimum with 1 here since $\mu_g$ is a bounded measure on $[0,1]^n$.
Given $(2.30)$ and the fact that first term in $(2.28)$ is bounded by $C|\lambda_0|^{-\epsilon}$, we conclude that
$$|T'(\lambda_0,...,\lambda_n)|  \leq C|\lambda_0|^{-\epsilon} + C |\lambda_0|^{-{1 \over v + 1}} \int_{|\lambda_0|^{-1}}^{\infty}t^{-{{1 \over v + 1} - 1}}\min(1, t^{\epsilon})  \,dt \eqno (2.31)$$
If $\epsilon < {1 \over v + 1}$, we use $\min(1, t^{\epsilon}) \leq t^{\epsilon}$ in $(2.31)$ and obtain that 
$|T'(\lambda_0,...,\lambda_n)|  \leq C|\lambda_0|^{-\epsilon}$, the desired estimate $(1.4)$. If $\epsilon = {1 \over v + 1}$
we gain an additional logarithmic factor. If $\epsilon > {1 \over v + 1}$, we use $1$ in the minimum if $t > 1$, and $t^{\epsilon}$ in the minimum if $t \leq 1$. In  this case we get
$$|T'(\lambda_0,...,\lambda_n)|  \leq C|\lambda_0|^{-\epsilon}  + C |\lambda_0|^{-{1 \over v + 1}} \int_{|\lambda_0|^{-1}}^1  t^{\epsilon - {1 \over v + 1} - 1}\, dt+ C |\lambda_0|^{-{1 \over v + 1}} \int_{1}^{\infty} t^{-{{1 \over v + 1} - 1}}  \,dt \eqno (2.32)$$
This is bounded by a constant times $|\lambda_0|^{-{1 \over v + 1}}$, again the desired estimate for part a) of Theorem 1.1.
Since the sharpness aspect  of part a) of Theorem 1.1 was proved at the end of section 1, we have completed the proof part a)  of Theorem 1.1.

\line{}

\noindent {\bf Proof of part b) of Theorem 1.1.} 

From now until the  end of the proof of part b) of Theorem 1.1 we use the notation $C$ to denote a constant that may depend on $f, g, $ and $p$ as in the statement of Theorem 1.1b). Note that part b) of Theorem 1.1 is immediate if $p = 1$ due to the integrability of $g$, so we assume that $p > 1$. Recall that $|T'(\lambda_0,...\lambda_n)|$ is the left-hand side of $(2.26)$. To bound
this for a given $(\lambda_0,...,\lambda_n)$, we proceed as follows. Let $l \geq 1$ be such that $|\lambda_l| \geq {1 \over n}
|(\lambda_1,...,\lambda_n)|$. Then $(2.26)$ implies that
$$|T'(\lambda_0,...\lambda_n)| \leq C\int_{[0,1]^n} \min\big(1,  \big|\sum_{\{m:\, a_m = a_l\}} \lambda_m x_m\big|^{-{1 \over v + 1}} \big) |g(x_1,...,x_n)|\, dx_1...\,dx_n\eqno (2.33)$$
By Holder's inequality, for any $1 <  p \leq \infty$ we therefore have
$$|T'(\lambda_0,...\lambda_n)| \leq C||g||_p \bigg(\int_{[0,1]^n} \min\big(1,  \big|\sum_{\{m:\, a_m = a_l\}} \lambda_m x_m\big|^{-{1 \over v + 1}}\big)^{p'}\, dx_1...\,dx_n\bigg)^{1 \over p'}\eqno (2.34)$$
Here ${1 \over p} + {1 \over p'} = 1$. For the fixed $l$ in $(2.34)$, let $\lambda''$ be the $n$-dimensional vector whose $m$th entry is $\lambda_m$ if $a_m 
= a_l$ and is zero otherwise. Note that $a_l$ was chosen so that $|\lambda''| \geq {1 \over n}|\lambda'|$. Let $v$ be the
unit vector in the direction of $\lambda''$. We perform the integration in $(2.34)$ first in the directions perpendicular to $v$ and then in the $v$ direction. We obtain
$$\int_{[0,1]^n} \big(\min\big(1,  \big|\sum_{\{m:\, a_m = a_l\}} \lambda_m x_m\big|^{-{1 \over v + 1}}\big)\big)^{p'}\, dx_1...\,dx_n \leq C\int_0^{\sqrt{n}}\big(\min\big(1, (|\lambda''|t)^{-{1 \over v + 1}}\big)\big)^{p'}\,dt \eqno (2.35)$$
$$\leq C|\lambda''|^{-1} + C|\lambda''|^{-{p' \over v + 1}}\int_{|\lambda''|^{-1}}^{\sqrt{n}} t^{-{p' \over v + 1}} \,dt \eqno (2.36)$$
If $p' < v + 1$, the above is bounded by $C|\lambda''|^{-{p' \over v + 1}}$. If $p' > v + 1$ it is bounded by  $C|\lambda''|^{-1}$, 
and if $p' = v + 1$ it is bounded by $C|\lambda''|^{-1}\ln|\lambda''|$. Put together, $(2.36)$ is bounded by 
$C|\lambda''|^{-\min({p' \over v + 1},1)}$ unless $p' = v + 1$ when one has an additional logarithmic factor. Since $|\lambda''|$ is
within a factor of $\sqrt{n}$ of $|\lambda'|$, we can replace $C|\lambda''|^{-\min({p' \over v + 1},1)}$ by $C|\lambda'|^{-\min({p' \over v + 1},1)}$ here. Inserting this back into $(2.34)$, if  $p' \neq  v + 1$ we see that
$$|T'(\lambda_0,...\lambda_n)| \leq C|\lambda'|^{-\min({1 \over v + 1},{1 \over p'} )}||g||_p \eqno (2.37)$$
If $p' = v + 1$ there is an additional factor of $|\ln |\lambda'||^{1 \over p'}$. This gives $(1.7a)$ for $p > 1$. (One can have $2 + |\lambda'|$ instead of $|\lambda'|$ in $(1.7a)$ simply because the absolute value of the integrand of 
$T'(\lambda_0,...\lambda_n)$ is integrable). This completes the proof of part b) of Theorem 1.1.

\line{}
\line{}

\noindent {\bf Proof of part c) of Theorem 1.1.} 

We repeat the argument from $(2.5)$ to $(2.12)$, setting $\lambda_0 = 0$.  At this point, we replace ${\bf M}$ by the $v$ by $v$ 
matrix ${\bf M'}$ given by deleting the first column and $v + 1$st row of ${\bf M}$. We then proceed as before. Since there is no longer a $\lambda_0F(v_1,...,v_n)r$ term in $(2.12)$, we
only need $v$ derivatives to make the argument work. Thus the remainder of the proof of part c) of the theorem proceeds exactly as part b), with $v + 1$ replaced by $v$ and we obtain
the estimate $(1.7b)$. This completes the proof of part c) of Theorem 1.1 and therefore the whole theorem.

\noindent {\bf Proof of Theorem 1.2.} 

\noindent When $a_i$ is a homogeneity appearing only once, $(2.33)$ becomes
$$|T'(\lambda_0,...\lambda_n)| \leq C\int_{[0,1]^n} \min\big(1,  |\lambda_i x_i|^{-{1 \over v + 1}} \big) |g(x_1,...,x_n)|\, dx_1...\,dx_n\eqno (2.38)$$
Equation $(2.38)$ is exactly $(2.27)$ with $\lambda_0$ replaced by $\lambda_i$ and $f(x_1,...,x_n)$ replaced by $x_i$. The
exact sequence of steps going from $(2.27)$ to $(2.32)$ therefore gives the estimate $(1.8a)$ in place of $(1.4)$ as desired. This gives
part a) of Theorem 1.2. As for part b), as in part c)  of the proof of Theorem 1.1, the absence of a $\lambda_0F(v_1,...,v_n)r$ term in $(2.12)$ when $\lambda_0 = 0$ means that one needs only $v$ derivatives in the uses of the Van der Corput lemma that led to $(2.38)$. As a result, 
we can improve the estimate $(2.38)$ to
$$|T'(0,\lambda_1...\lambda_n)| \leq C\int_{[0,1]^n} \min\big(1,  |\lambda_i x_i|^{-{1 \over v }} \big) |g(x_1,...,x_n)|\, dx_1...\,dx_n\eqno (2.39)$$
This time the steps from $(2.27)$ to $(2.32)$ give $(1.8b)$, which gives us part b) of Theorem 1.2. Since we dealt with the sharpness statement of part c) at the end of the last section, we are done
with the proof of Theorem 1.2.

\noindent {\bf Proof of Theorem 1.4.} 

Let $\psi(x_1,...,x_n)$ be a smooth nonnegative function supported on the unit ball and equal to 1 on a neighborhood of the
origin, and let $\psi_0(x_1,...,x_n) = \psi(2^{-{1 \over a_1}}x_1,...2^{-{1 \over a_n}}x_n) - \psi(x_1,...,x_n)$.
We rewrite the integral $(1.3)$ defining $T(\lambda_0,...,\lambda_n)$ as 
$$\int_{\R^n} e^{i\lambda_0 f(x_1,...,x_n) + i\lambda_1 x_1 + ... + i\lambda_n x_n}\phi(x_1,...,x_n)\psi(x_1,...,x_n)
dx $$
$$+ \sum_{j = 0}^K \int_{\R^n} e^{i\lambda_0 f(x_1,...,x_n) + i\lambda_1 x_1 + ... + i\lambda_n x_n}\phi(x_1,...,x_n)\psi_0(2^{-{j \over a_1}} x_1,...,2^{-{j \over a_n}}x_n)\, dx_1\,...\,dx_n \eqno (2.40) $$
Here $K$ denotes a constant depending on the size of the support of $\phi(x_1,...,x_n)$ and our estimates will not depend on $K$.
Note that the first term of $(2.40)$  satisfies the conditions of Corollary 1.3 with $g(x) = 1$, so the term is bounded in absolute
value by $C(1 + |\lambda|)^{-\epsilon}$ as needed. Hence we may devote our attention to bounding the sum in $(2.40)$.
We change variables in each term of the sum $(2.40)$ to $(y_1,...,y_n) = (2^{-{j \over a_1}} x_1,...,2^{-{j \over a_n}}x_n)$, and using the mixed  homogeneity of $f$ the sum becomes
$$\sum_{j = 0}^K\int_{\R^n} e^{i2^j \lambda_0 f(y_1,...,y_n) + i\lambda_1 2^{{j \over a_1}} y_1 + ... + i\lambda_n 2^{{j \over a_n}} y_n}(2^{j\sum_{i = 1}^n{1 \over a_i}})\phi(2^{j \over a_1}y_1,...,2^{j \over a_n}y_n)$$
$$\times \psi_0( y_1,..., y_n) \, dy_1\,...\,dy_n \eqno (2.41) $$
Note that due to the $\psi_0( y_1,..., y_n)$ factor the integrand in $(2.41)$ is supported on an annulus contained in the unit ball
but not intersecting a smaller ball centered at the origin.  We now apply 
Corollary 1.3 to each term in $(2.41)$ and add the resulting estimates. In order to do this, we determine
for each $j$th term a number $M_j$ for which the function $p_j(y_1,...,y_n) = (2^{j\sum_{i = 1}^n{1 \over a_i}})\phi(2^{j \over a_1}y_1,...,2^{j \over a_n}y_n)  \psi_0( y_1,..., y_n)$ satisfies $(1.2)$ with $g(x) = M_j$. As a result, by Corollary 1.3 the corresponding term of the sum in $(2.41)$ will be bounded in absolute value by $CM_j 2^{-j\epsilon}$. We will see  that $\sum_j M_j2^{-j\epsilon}$ is uniformly bounded 
when $(1.11)$ holds, and  therefore  Theorem 1.4 follows.

\noindent A direct computation using $(1.11)$ reveals that
$$|p_j(y_1,...,y_n)|, \,\,\,\,\,\,\,\sum_{i=1}^n \big|y_j {\partial p_j \over \partial y}(y_1,...,y_n)\big| \,\,\,\leq\,\,\, CA(2^{j\sum_{i = 1}^n{1 \over a_i}})(2^{-jk}) (\sum_{i=1}^n |y_i|^{a_i})^{-k}\eqno (2.42)$$
The constant $C$ here just depends on the function $\psi$. Since $(y_1,...,y_n)$ is in the annulus centered at the origin, the $(\sum_{i=1}^n |y_i|^{a_i})^{-k}$ factor is bounded, so we have
$$|p_j(y_1,...,y_n)|, \,\,\,\,\,\,\,\sum_{i=1}^n \big|y_j {\partial p_j \over \partial y}(y_1,...,y_n)\big| \,\,\,\leq\,\,\, C'A2^{j\big(\sum_{i = 1}^n{1 \over a_i}-k\big)} \eqno (2.43)$$
Hence $p_j(x)$ satisfies $(1.2)$ with $g(x) = C''A2^{j(\sum_{i = 1}^n{1 \over a_i}-k)}$. So as long as $k > \sum_{i = 1}^n{1 \over a_i }-\epsilon$, which we are assuming, then for $\zeta = k - \sum_{i = 1}^n{1 \over a_i } + \epsilon > 0$, Theorem 1.1 a) gives
 $$M_j 2^{-j\epsilon} \leq C'''2^{j\big(\sum_{i = 1}^n{1 \over a_i}-k\big)}\times 2^{-j\epsilon} = C2^{-j\zeta} \eqno (2.44)$$
 Thus $\sum_j M_j 2^{-j\epsilon}$ is uniformly bounded as needed. This completes the proof of Theorem 1.4.

\line {}

\noindent {\heading References.}

\line{}

\noindent [BNW] J. Bruna, A. Nagel, and S. Wainger, {\it Convex hypersurfaces and Fourier transforms},
Ann. of Math. (2) {\bf 127} no. 2, (1988), 333--365.  \parskip = 4pt\baselineskip = 4pt

\noindent [ChKNo] K. Cho, J. Kamimoto, T Nose, {\it Asymptotic analysis of oscillatory integrals via the Newton polyhedra of the phase and the amplitude}, J. Math. Soc. Japan 
{\bf 65} (2013), no. 2, 521-562. 

\noindent [CDMaM] M. Cowling, S. Disney, G Mauceri, and D. M\"uller {\it Damping oscillatory integrals}, 
Invent. Math. {\bf 101}  (1990),  no. 2, 237-260.

\noindent [CMa] M. Cowling, G. Mauceri, {\it Oscillatory integrals and Fourier transforms of surface carried measures},
Trans. Amer. Math. Soc.  {\bf 304} (1987), no. 1, 53-68.

\noindent [Gr] M. Greenblatt, {\it The asymptotic behavior of degenerate oscillatory integrals in two dimensions}, J. Funct. Anal. {\bf 257} (2009), no. 6 , 1759-1798. 

\noindent [GPT] A. Greenleaf, M. Pramanik, W. Tang, {\it Oscillatory integrals with 
homogeneous polynomial phase in several variables},  J. Funct. Anal. {\bf 244}  (2007),  no. 2, 444--487. 

\noindent [ISa] A. Iosevich, E. Sawyer, {\it Maximal averages over surfaces},  Adv. Math. {\bf 132} 
(1997), no. 1, 46--119.

\noindent [L] W. Littman, {\it Fourier transforms of surface-carried measures and differentiability of surface averages.}
Bull. Amer. Math. Soc. {\bf 69} (1963) 766-770. 

\noindent [NSeW] A. Nagel, A. Seeger, and S. Wainger, {\it Averages over convex hypersurfaces},
Amer. J. Math. {\bf 115} (1993), no. 4, 903-927.

\noindent [PY] M. Pramanik, C.W. Yang, {\it Decay estimates for weighted oscillatory integrals in ${\bf R}^2$},
Indiana Univ. Math. J., {\bf 53}  (2004), 613-645.

\noindent [RS] F. Ricci, E.M Stein  {\it Harmonic analysis on nilpotent groups and singular integrals. I. Oscillatory integrals.}
 J. Funct. Anal. {\bf 73} (1987), no. 1, 179–194.

\noindent [S] E. M.  Stein, {\it Harmonic analysis; real-variable methods, orthogonality, and oscillatory 
integrals}, Princeton Mathematics Series Vol. 43, Princeton University Press, Princeton, NJ, 1993.

\noindent [Sc] H. Schulz, {\it Convex hypersurfaces of finite type and the asymptotics of their Fourier 
transforms}, Indiana Univ. Math. J. {\bf 40} (1991), no. 4, 1267--1275. 

\noindent [SoS] C. Sogge and E. Stein, {\it Averages of functions over hypersurfaces in $R^n$}, Invent.
Math. {\bf 82} (1985), no. 3, 543--556.

\noindent [V] A. N. Varchenko, {\it Newton polyhedra and estimation of oscillating integrals}, Funct. Anal. Appl., 10 (1976), 175-196.

\line{}
\line{}

\noindent Department of Mathematics, Statistics, and Computer Science \hfill \break
\noindent University of Illinois at Chicago \hfill \break
\noindent 322 Science and Engineering Offices \hfill \break
\noindent 851 S. Morgan Street \hfill \break
\noindent Chicago, IL 60607-7045 \hfill \break
\noindent greenbla@uic.edu \hfill\break

\end